\documentclass[3p,number,review]{elsarticle}
\usepackage[latin9]{inputenc}
\usepackage{a4}
\usepackage{amsfonts}
\usepackage{amssymb}
\usepackage{amsmath}
\usepackage{amsthm}
\usepackage{nomencl}
\usepackage[colorlinks,linkcolor=blue,citecolor=green]{hyperref}
\begin{document}
\pagenumbering{arabic}
%
%
\newtheorem{thm}{Theorem}
\newtheorem{lem}[thm]{Lemma}
\newproof{pf}{Proof}
\theoremstyle{definition}
\newtheorem{defi}{Definition}
\title{On the number of cyclic transitive subgroups of a permutation group}
\author{Joachim K\"onig}
\address{Universit\"at W\"urzburg, Am Hubland, 97074 W\"urzburg, Germany}
\ead{joachim.koenig@mathematik.uni-wuerzburg.de}

\begin{abstract}
We prove an upper bound for the number of cyclic transitive subgroups in a finite permutation group (Note the update on p.\pageref{update}!) and clarify the structure of the groups for which this bound becomes sharp. 
We also give an application in the theory of number fields.\\
\end{abstract}
\begin{keyword}
Permutation groups; Galois groups
\end{keyword}
\maketitle
\section{Cyclic transitive subgroups in finite permutation groups}
Due to the classification of finite simple groups, the finite primitive permutation groups containing a cyclic transitive subgroup are known. The result, obtained by Feit (\cite{F}) and Jones (\cite{J}),
is the following:
\begin{thm}
\label{h1}
Let $G\leq S_n$ be a primitive group containing a cyclic transitive subgroup. Then one of the following holds:
\begin{enumerate}
\item $n=p \in \mathbb{P}$, $C_p \leq G \leq AGL_1(p)$.
\item $G = A_n$ or $S_n$.
\item $n = \frac{q^d-1}{q-1}$ with $d \geq 2$ and $q$ a prime power, $PGL_d(q) \leq G \leq P\Gamma L_d(q)$.
\item $n = 11$, $G = PSL_2(11)$ in its action on 11 points.
\item $G = M_{11}$ or $M_{23}$ in the natural action.
\end{enumerate}
\end{thm}
From this classification, one easily obtains a bound for the number of cyclic transitive subgroups in primitive permutation groups $G$ of degree $n$ (cf.\ the beginning of the next section):
This number is never larger than $|G|/n$, and if $G$ is not cyclic, it is even strictly smaller.\\
In the following, I will prove a generalization of this statement for arbitrary finite transitive permutation groups.
Apart from being of purely group-theoretic interest, this result has number-theoretic applications, see Section \ref{Appl}.
\begin{thm}[Main Theorem]
\label{haupt}
Let $G \leq S_n$ be a transitive permutation group on $n$ points. Then the number of cyclic transitive subgroups of $G$ is at most $\frac{|G|}{n}$. (Equivalently, the number of conjugacy classes of $n$-cycles in $G$ is at most $\varphi(n)$.)\\
Furthermore, if equality holds, then $G$ embeds into an iterated wreath product $(AGL_1(p_1) \wr AGL_1(p_2)) \wr ...\wr AGL_1(p_{m-1})) \wr C_{p_m}$, with $p_1,...p_m \in \mathbb{P}$ and $\prod p_i = n$. In particular, $G$ is then solvable.
\end{thm}
Before proceeding to the proof, we need an auxiliary lemma:
\begin{lem}
\label{h3}
Let $H$ be a primitive group of degree $n$ containing an $n$-cycle $\tau$. Then there is no proper cyclic overgroup of $\langle\tau\rangle$ in $Aut(H)$.
\end{lem}
\begin{pf}
The relevant cases are only the ones where $Aut(H)$ is not contained in $S_n$ (i.e. for our purposes only the cases $PGL_d(q) \leq H \leq P\Gamma L_d(q)$, $d\geq 3$; $A_6 \leq H \leq S_6$, 
and $H=PSL_2(11)$ acting on 11 points). Here the group of permutation automorphisms has index 2 in $Aut(H)$. So it suffices to show that in the transitive action of $Aut(H)$ on $2n$ points, 
there cannot be a $2n$-cycle. We prove this only for the only infinite series, and here we can assume $H = PGL_d(q)$, with $d \geq 3$. So let $\gamma$ be the graph automorphism of $PGL_d(q)$. 
If for some $x \in PGL_d(q)$, $x\gamma$ were cyclic transitive, then $x^\gamma$ would be a Singer cycle of $PGL_d(q)$. Progressing to inverse images in $GL_d(q)$, this means that for some 
$A \in GL_d(q)$, $A\cdot (A^t)^{-1} =:X$ is a pre-image of a Singer cycle in $PGL_d(q)$, i.e. the $k$-th power of a $GL_d(q)$-Singer cycle, $k \leq q-1$. 
Here $A^t$ denotes the transpose of $A$. However, as $X^t = A^{-1}A^t$ and $X^{-1} = A^t A^{-1}$ are conjugate in $GL_d(q)$ (and any element of $GL_d(q)$ is conjugate to its transpose), 
$X$ is conjugate to $X^{-1}$. But $X = Y^k$, $k \leq q-1$, for a Singer cycle $Y$, and the latter is conjugate in $GL_d(q)$ only to its $q^j$-th powers. 
As $q^{d-1} \cdot k \leq q^{d-1} \cdot (q-1) < q^d-1 - (q-1)$, $X$ cannot be conjugate to its inverse, and the assertion follows.
\end{pf}
\section{Proof of the Main Theorem}
First we show that the assertion of the theorem holds for all primitive groups (and the inequality becomes sharp only for the groups $C_p$). 
Using the list in theorem \ref{h1}, this is clear except for the case $PGL_d(q) \leq G \leq P\Gamma L_d(q)$. Here, with the one exception $G = P\Gamma L_2(8)$, all $n$-cycles are Singer cycles of $PGL_d(q)$ (cf. \cite{J}, Theorem 1).\\ In $P\Gamma L_2(8)$ there are $3 < \varphi(9) = 6$ conjugacy classes of $9$-cycles. In the other cases it is known that all Singer cyclic subgroups are conjugate, so the number of classes of $n$-cycles is $\frac{N_{S_n}(S)}{N_G(S)}$, $S$ being a Singer cycle. But $|N_G(S)| \ge |N_{PGL_d(q)}(S)| = n \cdot d > n$, see e.g. \cite{H}, $\mathsection$ 7, 7.3.\\
\\
So now let $G$ act imprimitively on $n$ points, with a partition of length $r$ and block length $s$ (i.e. $n=rs$), with $s > 1$ being minimal.
Then the homomorphic image $H$ of the block kernel $K$ in its action on a single block is transitive (as it contains an $s$-cycle), and it is a normal subgroup of a primitive group (namely the block stabilizer in the corresponding action), so by the list in theorem \ref{h1} $H$ is also primitive and contained in that list.\\
Also by induction, the fraction of $r$-cycles in the factor group $G/K$ cannot be higher than $\frac{\varphi(r)}{r}$.\\
Let $\tau$ be an $n$-cycle in $G$.\\ We distinguish between several cases:
%
\subsection{The case $s=p \in \mathbb{P}$, $C_p \leq H \leq AGL_1(p)$}
By the multiplicative properties of the Euler Phi-function, we can assume w.l.o.g that $p$ does not divide $r$, and even then it suffices to find a single element $x \in K$ such that $x\tau$ is not an $n$-cycle.\\
But $\tau^r \in K$, so the elements $\tau, \tau^{r+1},...,\tau^{(p-1)r + 1}$ are all in $K\tau$. Now if all of them were $n$-cycles, the exponents would all be coprime to $p$.
But this is obviously impossible, as there is $m \in \{1,...,p-1\}$ such that $mr=-1 \mod p$. 
%
\subsection{The case $S \le H \le Aut(S)$ with a non-abelian simple group $S$}
 Again, we first choose a fixed $n$-cycle $\tau$. By renumbering the elements in the different blocks, we can assume that $\tau$ acts as an $s$-cycle $\pi$ on the last block, followed by a cyclic permutation $(1...r)$ of the blocks. We use this setting only to find out more about the structure of $K$.\\
With $S:=soc(H)$, $K$ has as a normal subgroup a direct product $N$ of copies of $S$.
It holds that $|N| = |S|^{r'}$, where $r'$ divides $r$.  More specifically, the action of $N$ on e.g. the blocks $1, r/r'+1,...$ is diagonalized in the following sense: 
if $(\sigma_1,...,\sigma_r) \in N$, then \[\sigma_{k+r/r'} = \begin{cases}\sigma_k^{\pi\rho}, \text{ if } k \le r < k+r/r'\\ \sigma_k^\rho, \text{ otherwise}\end{cases}\] with an automorphism 
$\rho \in Aut(S)$.\footnote{Of course, the index in the above formula is always to be regarded modulo $r$} 
Here $\rho$ is independent of $k$ because of the conjugation action of $\langle\tau\rangle$ on $K$. 
\\(
To see this, consider the extension $N\tau$. The group $N$ has only $r'$ normal subgroups of order $|S|^{r'-1}$, and these are all kernels of suitable projection homomorphisms of $N$ onto a certain block. The existence of a cyclic transitive subgroup in the action on the blocks assures that all kernels appear with the same frequency, so for $r/r'$ blocks each it holds that the image of an element of $N$ on a given one of those blocks is 1, if and only if the images on all the other ones are 1. So the images on these blocks differ only by an automorphism of $S$, i.e. $\sigma_{1+r/r'} = \sigma_1^\rho$, and the rest of the assertion follows immediately by the action of the cycle $\tau$, as indicated above.
)\\
\\
It is now easy to see that the same diagonalization property holds for all $(\sigma_1, ..., \sigma_r) \in K$. Namely, for all $(\widehat{\sigma}_1, ..., \widehat{\sigma}_r) \in N$, it holds that $(\widehat{\sigma}_1^{\sigma_1}, ..., \widehat{\sigma}_r^{\sigma_r}) \in N$, therefore e.g. $\widehat{\sigma}_1^{\sigma_1\rho} = \widehat{\sigma}_{1+r/r'}^{\sigma_{1+r/r'}}$, and consequently $\sigma_1\rho\sigma_{1+r/r'}^{-1} = \rho$, i.e. $\sigma_{1+r/r'} = \sigma_1^\rho$.\\
\\
In particular, $\pi\rho^{r'} = 1$; but as $\pi$ is an $s$-cycle, this means that $r'$ (= the number of blocks belonging to one diagonalized component) and $s$ are coprime.\\
Now we will use slightly different setups, depending on whether $r$ and $s$ are coprime or not.
\\ 
\subsubsection{The case $(r,s) = 1$}
Here the product $N.\tau^s$ splits, and so we can assume that $\tau^s$ is simply a cyclic permutation $(1...r)$ of the blocks.\\
Then there is $\rho \in Aut(S)$ such that $\sigma_{k+r/r'} = \sigma_1^\rho$ for all $(\sigma_1,...\sigma_r) \in K$, and  {\it{all}} $k \in \{1,...,r\}$.
\\
\\
Now, if $(\sigma_1,...\sigma_r) \cdot \tau^s \in K\tau^s = K\tau$ is an $n$-cycle, the projection of its $r$-th power to the first component is 
\[(\sigma_1 \cdots \sigma_{r/r'} \cdot \rho^{-1})^{r'}\]
and this element is then an $s$-cycle in $H$, which firstly, by lemma \ref{h3}, means that $\rho^{-1} \in H$, and secondly is equivalent to $\underbrace{\sigma_1\cdots\sigma_{r/r'}}_{=:x} \cdot \rho^{-1}$ being an $s$-cycle, as $(r,s)=1$.\\
\\
Now how often does a given $x$ occur as a product $\sigma_1\cdots\sigma_{r/r'}$, if we run through all of $K$?
To see this, we further split up $\sigma_i:=\alpha_i\cdot \beta_i$ into an $S$-part $\alpha_i$ and an outer automorphism part $\beta_i$. So
\[\sigma_1\cdots\sigma_{r/r'} \cdot \rho^{-1}= \alpha_1 \cdot \alpha_2^{\beta_1^{-1}} \cdots \alpha_{r/r'}^{...} \cdot (\beta_1 \cdots \beta_{r/r'}) \cdot \rho^{-1}\]

To count $n$-cycles in $K\tau$, we can now fix all other variables in the product while running through all of $S$ with $\alpha_1$. Let $c_H$ denote the number of $s$-cycles in $H$.\\
In cases where $S=H$ is a simple group, it is obvious that the number of $n$-cycles in $K\tau$ is bounded by $\frac{|K|}{|H|} \cdot c_H$. Therefore the total number of $n$-cycles in $G$ cannot exceed $\frac{|K|}{|H|} \cdot c_H$ times the number of $r$-cycles in $G/K$, and by induction this last number does not exceed $[G:K] \cdot \frac{\varphi(r)}{r}$. So there are at most $|G| \cdot \frac{c_H}{H} \cdot \frac{\varphi(r)}{r} < |G| \cdot \frac{\varphi(s) \varphi(r)}{sr} \le |G| \cdot \frac{\varphi(n)}{n}$ $n$-cycles in $G$.\\
\\
In case $S=A_s$, $H=S_s$, the $s$-cycles of $H$ lie either all inside or all outside $S$ (depending on whether $s$ is odd or even). In the first case, at most $\frac{2|S|}{s}$ choices for $\alpha_1$ can lead to $\sigma_1 \cdots \sigma_{r/r'} \cdot \rho^{-1}$ being an $s$-cycle (again, all other variables being fixed). As $2 < \varphi(s)$ ($s\ge 5$ is odd!), the number of $n$-cycles in $G$ is strictly less than $|G| \cdot \frac{\varphi(n)}{n}$, just as in the previous case. For $s$ even, note also that $\rho^{r'} = 1$, and as $r'$ is odd, $\rho$ lies inside $A_s$. So the outer automorphism part of $\sigma_1\cdots\sigma_{r/r'} \cdot \rho^{-1}$ is just the product of the $\beta_i$, which yields a homomorphism into $C_2$. Therefore at most half of the possible choices for the $\beta_i$ can lead to elements outside $A_s$, and for those, again at most $\frac{2|S|}{s}$ choices for $\alpha_1$ can lead to $\sigma_1 \cdots \sigma_{r/r'} \cdot \rho^{-1}$ being an $s$-cycle, i.e. the number of $n$-cycles in $G$ is at most $\frac{|K|}{s} \cdot [G:K] \cdot \frac{\varphi(r)}{r} < |G| \cdot \frac{\varphi(n)}{n}$.\\
\\
Finally, let $S=PSL_d(q)$. Then no $s$-cycles of $H$ lie inside a proper normal subgroup of $PGL_d(q)$. But also none of them lie outside of $PGL_d(q)$, except in the case $s=9$, $H=P\Gamma L_2(8)$, as mentioned earlier.\\ 
We leave out the one exceptional case for now.
Note that $\rho$, because of its order being coprime to $s$, has to lie in $P\Sigma L_d(q)$, i.e. in order for $x\rho^{-1}$ to be an $s$-cycle, the product of the $\beta_i$ has to lie in $\rho PGL_d(q)/PSL_d(q)$.\\

Again, if everything else is fixed, at most $\underbrace{\frac{\varphi(s)}{d \cdot s} \cdot |S| \cdot (d,q-1)}_{\text{number of $s$-cycles in $PGL_d(q)$}}$ choices for $\alpha_1$ can lead to an $s$-cycle.\\
Now if {\it all} the choices for the $\beta_i$ would give rise to $s$-cycles (for fitting $\alpha_i$), we would get at most $\frac{|K|}{|S|} \cdot \frac{\varphi(s)}{d \cdot s} \cdot |S| \cdot (d,q-1) = |K| \cdot \frac{\varphi(s)}{s} \cdot \frac{(d,q-1)}{d}$ $n$-cycles in $K\tau$.\\
\\
But if $\rho$ lies outside $PSL_d(q)$, then, using that $H/S$ is metacyclic, 
we see that not all products $\beta_1 \cdots \beta_{r/r'}$ can lie in the $\rho$-coset, whereas if $\rho \in PSL_d(q)$, out of all those products $\beta_1 \cdots \beta_{r/r'}$ inside $PGL_d(q)/PSL_d(q)$, at most a fraction of $\frac{\varphi((d,q-1))}{(d,q-1)}$ has $PGL_d(q)/PSL_d(q)$-part of maximal order and can therefore give rise to an $s$-cycle (when combined with the right element inside $S$).\\
In any case, the number of $n$-cycles in $K\tau$ is therefore strictly smaller than $|K| \cdot \frac{\varphi(s)}{s}$, which proves the assertion for this case.
\\ \\
In the special case $H=P\Gamma L_2(8)$, there is one $H$-class of $n$-cycles in each coset of $S$. So again iteration over $\alpha_1 \in S$ yields an upper bound 
$|K| \cdot \frac{|H|}{9 \cdot |S|} < |K| \cdot \frac{\varphi(9)}{9}$ for the number of $n$-cycles in $K\tau$.\\
\subsubsection{The case $(r,s) \ne 1$}
Then the extension $K\tau$ need not split. We fix an $n$-cycle $\tau$ and again count $n$-cycles while running through $K\tau$. By renumbering (i.e. permuting) the elements in the blocks, we can assume (as $(r',s) = 1$!) that $\tau$ acts on {\it all} the blocks $k \cdot r/r'$ ($k=1,...,r'$) as the same $s$-cycle $\pi$, followed by permuting the blocks cyclically.\\
Then, if $\rho \in Aut(S)$ such that $\sigma_{1+r/r'} = \sigma_1^\rho$ for all $(\sigma_1,...\sigma_r) \in K$, the action of $\tau$ shows that for all $k \in \{1,...,r\}$: $\sigma_{k+r/r'} = \sigma_k^{\rho^{(\pi^j)}}$ for suitable $j \in \mathbb{N}_0$.\\
But also $(\sigma_1,...,\sigma_r)^{\tau^r} = (\sigma_1^{\pi^{r'}},...,\sigma_r^{\pi^{r'}})$, so $\rho$ commutes with $\pi^{r'}$, and therefore also with $\pi$.\\
This, however, means that $\sigma_{k+r/r'} = \sigma_k^\rho$ for {\it all} $k \in \{1,...,r\}$.
\\
Now if $(\sigma_1,...\sigma_r) \cdot \tau$ is an $n$-cycle, the projection of its $r$-th power to the first component is 
\[(\sigma_1 \cdots \sigma_{r/r'} \cdot \pi\rho^{-1})^{r'}\]
This element is now required to be an $s$-cycle (which is equivalent to $\sigma_1 \cdots \sigma_{r/r'} \cdot \pi\rho^{-1}$ being an $s$-cycle, as $(r',s)=1$). We can split this product up into $S$-part and outer automorphism part, just as above. The only difference is that this time the outer automorphism part of $\pi\rho^{-1}$ can be arbitrary, which increases the bound on number of $n$-cycles in $K\tau$ slightly in some cases, namely to $\frac{|K|}{|S|} \cdot c_H$ in case $S=A_s$ with even $s$ (because it need no longer be true that half of the products obtained in the above way must lie inside $A_s$!), and to 
\[\frac{|K|}{|S|} \cdot c_H \leq \frac{|K|}{|S|}  \cdot \frac{\varphi(s)}{d \cdot s} \cdot |S| \cdot (d,q-1) = |K| \cdot  \frac{\varphi(s)}{s} \cdot \frac{(d,q-1)}{d}\] in the case $S=PSL_d(q)$ (because analogously, it need no longer hold that a certain fraction of products must lie inside $PSL_d(q)$). \\
So in all cases, the number of $n$-cycles in $K\tau$ is not larger than $|K|\cdot \frac{\varphi(s)}{s}$, which is enough for our purposes, as now $\varphi(s)\varphi(r) < \varphi(sr) = \varphi(n)$.
\subsection{The case $H = S_4$}
Here, $K$ has as a normal subgroup an elementary-abelian 2-group $N$, and for any $n$-cycle $\tau$, $\tau^{n/2}$ lies in $N$. As before, in the case $(r,4) = 1$, we get as a necessary condition for $\sigma\tau^r \in K\tau$ to be an $n$-cycle that $\sigma_1 \cdots \sigma_r$ must be a 4-cycle. Again, we split $\sigma_i = \alpha_i\cdot\beta_i$, with $\alpha_i \in Soc(S_4)$, and $\beta_i$ an outer automorphism of the socle. As $N$ is elementary-abelian, for fixed $\beta_i$'s the product map $N$ yields a homomorphism $(\alpha_1,...,\alpha_r) \mapsto \prod{\alpha_i}$ into $C_2 \times C_2$. However for fixed $\beta_i$'s, at most half of the possible choices in the socle can have a product leading to a 4-cycle (namely those with product one of two certain involutions of $C_2 \times C_2$). In the same way, at most half of the possible choices of the $\beta_i$ can lead to a product outside $Alt(3)$, which is also necessary for a 4-cycle. So the number of $n$-cycles in $K\tau$ for a given $\tau$ is at most $\frac{|K|}{4}$, whereas $\frac{|K|}{2}$ would be required to reach the critical value.\\
In case $(r,4) >1$ it suffices to find a single element in $K$ that does not give rise to an $n$-cycle in $K\tau$.  Letting $\tau$ act as a product of $\pi = (1 2 3 4)$ (on the last block), and a cyclic permutation $(1...r)$ of the blocks, we get the condition that $\sigma_1 \cdots \sigma_r \pi$ must be a 4-cycle for all $(\sigma_1,...\sigma_r) \in K$. Again writing $\sigma_i = \alpha_i\beta_i$, we have
\[\sigma_1\cdots\sigma_r \cdot \pi = \underbrace{\alpha_1\cdot \alpha_2^{\beta_1^{-1}} \cdots \alpha_r^{...}}_{=:p_\alpha} \cdot \underbrace{(\beta_1 \cdots \beta_r)}_{=:p_\beta} \cdot \pi\]
Let us look only at $\sigma \in K \cap Alt(4)^r$. Then $p_\beta$ yields a homomorphism into $C_3$, and for fixed $\beta_i$, $p_\alpha$ yields a homomorphism $Soc(K) \to C_2\times C_2$ (the $\alpha$'s and $\beta$'s can of course be chosen independently, as the extension splits).
Now $p_\beta$ cannot always be 1; for $\sigma \in K \cap Alt(4)^r$ with $\beta_r = (1 2 3)$ and $p_\beta = 1$, we have $p_\beta(x^\tau) \ne 1$. 
Thus there exists $\sigma$ with $p_\alpha(\sigma)=1$, $p_\beta(\sigma) = (1 3 2)$ - and then the product above is not a 4-cycle any more. This completes the proof.\\\\
{\bf Remarks:}
\begin{enumerate}
\item As the bound $\varphi(n)$ on the number of conjugacy classes of $n$-cycles seems so natural, it would be nice to have a more direct proof for this statement. 
However, intuitive approaches seem to fail, e.g. there is no direct relation between the size of an $n$-cycle normalizer and the number of $n$-cycle classes in $G$. 
On the contrary, there are even examples where the bound becomes sharp, but still the whole normalizer $AGL_1(n)$ is contained in $G$, at least for some of the $n$-cycles in $G$. 
An infinite series of examples for this is provided by groups of the form $((C_{3^k} \times C_{3^k}).Aut(C_{3^k})).2$, of degree $2\cdot3^k$, namely the subgroups of $AGL_1(3^k) \wr C_2$, 
where the $Aut(C_{3^k})$-part in the block kernel has been diagonalized. (The number of $n$-cycles in these groups can easily be verified similarly as above).
\item In the maximal case as well as in the minimal one (where $G$ has just one conjugacy class of $n$-cycles), it is true that the number of cyclic transitive subgroups of $G$ divides $|G|$. This is however not true in general; the smallest counter-example is $C_3\wr C_3$, of order 81, which has 6 cyclic transitive subgroups.
\\
\end{enumerate}

\section{A number theoretic application}
\label{Appl}
Let $K$ be a number field, $f \in K[X]$ be irreducible of degree $n$, and $G = Gal(f|K)$, $x$ a root of $f$ over $K$.
By Chebotarev's density theorem, the fraction of $n$-cycles in $G$ asymptotically equals the fraction of prime ideals $p$ of $K$ that remain prime in $K(x)$, 
which in turn is equivalent to $f$ remaining irreducible modulo $p$. So theorem \ref{haupt} says that the density of such primes is at most $\frac{\varphi(n)}{n}$. 
Furthermore, if $G$ is unknown, and out of a large number $N$ of randomly chosen primes, $f$ remains irreducible modulo $p$ for approximately $N \cdot \frac{\varphi(n)}{n}$ of them,
then $G$ is very likely solvable.\\
\\
Using some well-known results about monodromy groups of polynomials (cf. e.g. \cite{M} or \cite{Fr} for an introduction), one can easily construct polynomials $f \in \mathbb{Q}[X]$ such that $Gal(f(X)-t|\mathbb{Q}(t)) = ((C_{3^k} \times C_{3^k}).Aut(C_{3^k})).2$, the groups mentioned above, e.g. $f = (X^{3^k}-1)^2$ 
(Note that $AGL_1(n) \le G$ is a necessary condition for $G$ to be the monodromy group of a polynomial over $\mathbb{Q}$, by a special case of Fried's "branch cycle argument".).\\
Hilbert's irreducibility theorem then yields infinitely many specializations $t_0 \in \mathbb{Q}$ such that the Galois group is preserved, and so the set of prime numbers $p \in \mathbb{P}$ 
for which $f(X)-t_0$ remains irreducible modulo $p$ has maximal possible density, in the sense of Theorem \ref{haupt}.\\
\\
{\bf Acknowledgement:}\\
I would like to thank Lior Bary-Soroker and Peter M\"uller for pointing this problem out to me.
\vspace{10mm} \ \\
{\bf Update:}
\label{update}\\
Thanks to Ilia Ponomarenko for informing me that the first part of the statement of Theorem \ref{haupt} has already been obtained by Mikhail Muzychuk in \cite{Mu}.

\end{document}